\documentclass[reqno]{amsart}
\usepackage{amsmath}
\usepackage{amssymb}
\usepackage{amsthm}
\usepackage{enumerate}
\usepackage[mathscr]{eucal}
\usepackage{graphicx}

\newtheorem{theorem}{Theorem}[section]

\newtheorem{proposition}[theorem]{Proposition}
\newtheorem{corollary}[theorem]{Corollary}
\newtheorem{J-S Theorem}[theorem]{J-S Theorem}
\newtheorem{definition}[theorem]{Definition}
\theoremstyle{definition}
\newtheorem{example}[theorem]{Example}
\newtheorem{remark}[theorem]{Remark}

\begin{document}
	\title[On UC Quasi-Metric Spaces]{Quasi-metric spaces on which real-valued continuous functions are uniformly continuous}
	\author{Om Dev Singh and Anubha Jindal}
	
	\thanks{Om Dev Singh: Department of Mathematics, Malaviya National Institute of Technology Jaipur, Jaipur, Rajasthan, India, email: om2947dev@gmail.com}
	\thanks{Anubha Jindal: Department of Mathematics, Malaviya National Institute of Technology Jaipur, Jaipur, Rajasthan, India, email: anubha.maths@mnit.ac.in}
	
	
	\maketitle
	
	\begin{abstract}
		
		The concept of a quasi-metric space arises by relaxing the requirement of the symmetry axiom in the definition of a metric. This small variation alters several structural properties possessed by a standard metric space. This article aims to investigate the notion of UC quasi-metric spaces in a systematic manner. A quasi-metric space $(X,d)$ is called a UC space if every real-valued continuous function on $(X,d)$ is uniformly continuous. In the context of metric spaces, UC spaces help in bridging the gap between compactness and completeness. These spaces also play an important role in the theory of hyperspaces of closed sets and fixed point theory. In this article, we present several characterizations of UC quasi-metric spaces and provide various examples of such spaces. At several instances, our proof techniques highlight key differences between UC quasi-metric spaces and their metric counterparts. 
\end{abstract}	
		

\section{Introduction}
In 1931, Wilson introduced the notion of a quasi-metric space (\cite{Wilson}).  A quasi-metric $d$ on a set $X$ is a function on $X\times X$ that satisfies all the metric axioms, but the symmetry condition $d(x,y) =d(y,x)$ may fail for $x,y \in X$.  This subtle modification in the definition of a metric substantially changes the whole theory, specially with respect to completeness and compactness. For instance, there are several non-equivalent notions of Cauchy sequences, and consequently of completeness in quasi-metric spaces. Most of these notions coincide in the presence of metric structure. 	
The monograph \cite{fletcherlindgren} by Fletcher and Lindgren contains all the necessary fundamental aspects of quasi-uniform and quasi-metric spaces. A more recent treatment is provided by Cobzaș \cite{Cobzas}, who adopts a modern perspective on quasi-metric spaces. In particular, his monograph offers a detailed study of asymmetric normed spaces. Quasi-metric spaces are also sometimes referred to as asymmetric metric spaces or oriented metric spaces \cite{Bodjanova,Zimmer,ILR}.

In the last two to three decades, significant progress has been made in the study of quasi-metric spaces, with particular emphasis on function spaces, hyperspaces, and fixed point theory (see, \cite{RomagueraPedro_fp, cow1, kunjiromaguera, RomagueraLopez,Seceleanetal_fp}) and references therein. However, the study of various structural properties of quasi-metric spaces corresponding to different classes of functions remains largely undeveloped. In the realm of metric spaces, a large literature has been devoted to develop such studies. In fact, in metric space theory, functions often serve as tools to capture and reflect the intrinsic geometry as well as topology of the space.  Due to the absence of the symmetry property, quasi-metric spaces exhibit subtleties that complicate such a study in this case. Consequently, several classical results need to be reformulated or reinterpreted when carried over to this asymmetric framework. This gap highlights the need and provides the motivation behind the present work. Recently, in \cite{ODS_AJ_Cauchy_regular}, the authors studied relationship between various properties of quasi-metric spaces through \textit{left K-Cauchy regular maps} defined between them.  It is easy to see that every uniformly continuous map between quasi-metric spaces is left K-Cauchy regular. 
	
It is well-known that every real-valued continuous function defined on a compact metric space is uniformly continuous. Notably, several such consequences of compactness can be established even when the underlying metric space satisfies less restrictive hypotheses. As a result, many classes of metric spaces that lie between complete and compact metric spaces have been defined and extensively studied in the literature \cite{beer1, kundulipsymanisha}. An important such class is the class of UC metric spaces. A metric space is called a \textit{UC space} if every real-valued continuous function defined on it is uniformly continuous. UC metric spaces are also known as Atsuji spaces and have been well-studied in the literature. In particular, various equivalent characterizations of such spaces can be found in \cite{atsuji1958, beer1985, hueber1981, kundumanisha, monteiro1951, mrowka1965, nagata1950, toader1978, waterhouse1965}. For a well-written survey article on UC spaces, we refer readers to Jain and Kundu \cite{Kundujain2006}. The UC property has been used to study the equivalence of Vietoris and Hausdorff metric topologies on the collection of all closed subsets of a metric space (\cite{beer1985, beer1993}). Since every UC metric space is complete, these spaces are also useful in fixed point theory (see, \cite{beer1986b}). For recent advancements and the relationship of UC spaces with several other spaces and metric properties, one can refer to the recent monographs by Beer \cite{beer1} and Kundu et al. \cite{kundulipsymanisha}. 
	
The quasi-metric spaces on which every continuous real-valued function is uniformly continuous have been implicitly studied in the literature. For example, in \cite{MarinRomaguera1996} Marín and Romaguera proved that if a quasi-metric space $(X,d)$ has the Lebesgue number property (see, Definition \ref{Lebesgue number property defi}), then every continuous function from $(X,d)$ to the space $\mathbb{R}$ of all real numbers equipped with the usual metric is uniformly continuous. Moreover, it is known that every continuous function from a compact quasi-metric space to the space of reals with the usual metric is uniformly continuous. However, to the authors' knowledge, no complete characterization is known for those quasi-metric spaces on which every real-valued continuous function is uniformly continuous.    
	
The aim of this paper is to explore systematically the UC property in the setting of quasi-metric spaces. It is important to mention that the techniques used in the investigations of UC metric spaces are not necessarily applicable when transitioning to a quasi-metric setting (attributed to the lack of symmetry in such spaces). For instance, in a metric space, if a sequence $(x_n)$ is parallel to a sequence $(y_n)$ and $(x_n)$ fails to possess a cluster point, then $(y_n)$ also cannot admit a cluster point. This correspondence, however, does not necessarily hold in quasi-metric spaces. Furthermore, metric spaces enjoy very strong separation properties (often implying results about continuity). In contrast, a quasi-metric space need not be even $T_1$. These structural differences suggest that extending the notion of UC spaces from metric spaces to quasi-metric spaces is both mathematically interesting and demanding.
	
We begin in Section 2 by providing the necessary fundamental definitions and preliminary results. The core of our work starts in Section 3, where we formally define the UC property for quasi-metric spaces and explore some fundamental properties of such spaces. Subsequently, in Section 4, we provide several important characterizations of UC quasi-metric spaces, particularly in terms of isolation functional, uniformly discrete sets, and the Lebesgue number property. In Section 5, we first give a sequential characterization of UC quasi-metric spaces in terms of \textit{pseudo left K-Cauchy sequences}, an adaptation of the notion of pseudo Cauchy sequence, in the setting of quasi-metric spaces. Then we use it to provide a Cantor intersection type condition that characterizes UC quasi-metric spaces. Finally, we extend a recent result of Beer and Garrido \cite{beerGarridoMerono2018} that applies UC spaces to study when a reciprocal of a (non-vanishing) uniformly continuous function is uniformly continuous.

Since quasi-metric spaces often lacks strong separation axioms, when proving our results, we typically have to assume specific topological conditions. Consequently, our proofs clearly reveal that our findings are dependent on both the specific quasi-metric used and the topology it induces. To support our results, we have given several examples and counterexamples throughout the paper. Our investigations offer several new insights into the theory of quasi-metric spaces.

\section{Preliminaries}
	In this section, we summarize some basic definitions and results related to quasi-metric spaces. For more details, readers may refer to \cite{Cobzas, Zimmer}.
	
	\begin{definition}\normalfont
		A \textit{quasi-metric} on a set $X$ is a function $d: X \times X \rightarrow \mathbb{R}$ satisfying the following properties for $x, y, z \in X$
		
		\begin{enumerate}
			\item  $d(x, y) \geq 0$, $d(x,x)=0$
			\item  $d(x, y) = 0$ = $d(y,x)$ $\implies$ $x = y$	
			\item  $d(x, y) \leq d(x, z)+ d(z, y)$.
		\end{enumerate}
	\end{definition}
	
	If $d$ is a quasi-metric, then the \textit{conjugate quasi-metric} $\bar{d}$ is defined as $\overline{d}(x,y)=d (y,x)$ for $ x,y\in \;X.$
	Moreover, $d^{s}(x,y)=\max\{d(x,y),\overline{d}(x,y)\}$ for $x, y \in X$ is a metric on $X$, this metric is called the \textit{associated metric} to the quasi-metric $d$.
	
	\begin{definition}\normalfont
		Suppose $(X,d)$ is a quasi-metric space. Then
		\begin{enumerate}
			\item $B^{+}(x,\epsilon)=\{y\in X:d(x,y)<\epsilon\}$ is called the \textit{forward open ball}; and
			\item $B^{-}(x,\epsilon)=\{y\in X:d(y,x)<\epsilon\}$ is called the \textit{backward open ball}.
		\end{enumerate}  
	\end{definition}
	The topology $\tau(d)$ (resp. $\tau(\overline{d})$) on $X$ generated by all forward balls (resp. backward balls) is called the \textit{forward topology} (resp. \textit{backward topology}).

	We now give some important examples of quasi-metric spaces.
	
	\begin{example}\label{Example upper and lower quasi metrics on R}
		The function $d:\mathbb{R}\times \mathbb{R} \to \mathbb{R}$  defined by 
		$d(x,y)=\begin{cases}
			y-x, \text{ if } y\geq x\\    
			0, \text{   if }  y<x
		\end{cases}$\\
		is a quasi-metric on $\mathbb{R}$. Here $B^{+}(x,\epsilon)=(-\infty,x+\epsilon)$ and $B^{-}(x,\epsilon)= (x-\epsilon,\infty)$.
		The quasi-metric $d$ is called the \textit{upper quasi-metric} and its conjugate is called the \textit{lower quasi-metric} on $\mathbb{R}$. Occasionally, we may use the notations $d_l$ and $d_u$ for the lower and upper quasi-metrics on $\mathbb{R}$, respectively. It is easy to see that the associated metric to upper and lower quasi-metrics is the usual metric on $\mathbb{R}$. \qed 
	\end{example}
	Throughout this paper, the set $\mathbb{R}$ of all real numbers equipped with the usual metric is denoted by $(\mathbb{R},|\cdot|)$.
	\begin{example}\label{Example Sorgenfrey quasi metric} 
		The function	$d :\mathbb{R} \times \mathbb{R} \rightarrow  \mathbb{R}$ defined by
		$ d(x,y)=\begin{cases}
			y-x, \text{ if } y\geq x\\    
			1, \text{   if }  y<x
		\end{cases}$\\
		is a quasi-metric on $\mathbb{R}$. The forward topology on $\mathbb{R}$ induced by $d$ is the well-known lower limit or Sorgenfrey topology. Consequently, we call this quasi-metric  on $\mathbb{R}$, the \textit{Sorgenfrey quasi-metric}. \qed
	\end{example}

	Note that every quasi-metric space is first countable as well as $T_0$, however they need not be $T_1$ (see, Example \ref{Example upper and lower quasi metrics on R}). A quasi-metric space $(X,d)$ is $T_1$ if and only if $d(x,y)>0,$ whenever $x\neq y$. Also, $\tau(d)$ is $T_1$ if and only if $\tau(\bar{d})$ is $T_1$ (\cite{Cobzas}).

	From Example \ref{Example Sorgenfrey quasi metric}, it follows that separability and second countability are not necessarily equivalent in quasi-metric spaces. 
	\begin{definition}\normalfont
		A sequence $(x_{n})$ in a quasi-metric space $(X,d)$ is said to be \textit{forward convergent} (resp. \textit{backward convergent}) if there exists $ x\in X$ such that $\lim_{n}d(x,x_{n})=0$ (resp. $\lim_{n}d(x_{n},x)=0$). If $(x_n)$ is forward convergent (resp. backward convergent) to $x$, we denote it by $x_{n}\xrightarrow{f}x$ (resp. $x_{n}\xrightarrow{b}x$).
	\end{definition}
	\begin{enumerate}
		\item In quasi-metric spaces forward convergence may not imply backward convergence. For instance in the Sorgenfrey quasi-metric space (Example \ref{Example Sorgenfrey quasi metric}), the sequence $(1/n)$ is forward convergent to $0$ but it is not backward convergent to $0$. 
		\item In a $T_1$ quasi-metric space $(X,d)$, if a sequence $(x_n)$ is forward convergent to $x\in X$ and backward convergent to $y\in X$, then $x=y$ (see, Lemma 3.1 in \cite{Zimmer}).
	\end{enumerate}	
	
	Quasi-metric spaces having the property forward convergence implies backward convergence are also referred to as \textit{small set symmetric} in literature \cite{kunzireilly}. 
		
	\begin{definition}\normalfont{
			Let $(X, d)$ be a quasi-metric space and $(x_n)$ be a sequence in $X$. A point $x \in X$ is called a cluster point of the sequence $(x_n)$ in $(X,d)$ if there exists a subsequence $(x_{n_k})$ of $(x_n)$ such that $x_{n_k} \xrightarrow{f} x$.}		
	\end{definition}
	
	\begin{definition}\normalfont Let $(X,d)$ be a quasi-metric space. A sequence $(x_n)$ in $(X,d)$ is called \textit{left K-Cauchy} (resp. \textit{right K-Cauchy}) if for every $\epsilon>0$ there exists  $n_{0} \in \mathbb{N}$ such that $\forall$ $ n_{0} \leq k \leq n$, we have $d(x_{k},x_{n})<\epsilon$ (resp. $d(x_{n},x_{k})<\epsilon$).
	\end{definition}
	
	\begin{remark}\label{conv_not Cauchy}
		In metric spaces, every convergent sequence is Cauchy, while in a quasi-metric space, a forward convergent sequence need not be left K-Cauchy. For instance, in the Sorgenfrey quasi-metric space (Example \ref{Example Sorgenfrey quasi metric}), the sequence $(1/n)$ is forward convergent to $0$ but it is not left K-Cauchy. However, it can be easily shown that if a left K-Cauchy sequence $(x_n)$ has a cluster point in $(X,d)$, then $(x_n)$ is forward convergent.
		
	\end{remark}
	\begin{definition}\normalfont
		A quasi-metric space $(X,d)$ is called \textit{left K-complete} if every left K-Cauchy sequence is forward convergent.
	\end{definition}
	
	The upper and lower quasi-metrics on $\mathbb{R}$ as defined in Example \ref{Example upper and lower quasi metrics on R} are left K-Complete. However, the Sorgenfrey quasi-metric space (Example \ref{Example Sorgenfrey quasi metric}) is not left K-complete as $(x_n)=(\frac{-1}{n})$ is left K-Cauchy but not forward convergent.

		For quasi-metric spaces $(X,d)$ and $(Y,\rho)$, there are four types of continuities for a function $f:X\to Y$.
		\begin{definition}\normalfont  A function $f: (X,d) \to (Y,\rho)$ between two quasi-metric spaces $(X,d)$ and $(Y,\rho)$ is called forward-forward continuous (ff-continuous) at $x_0 \in X$, if for every $\epsilon >0$ there exists $\delta >0$ such that $\rho(f(x_0), f(y)) < \epsilon$ whenever $d(x_0,y) < \delta$.
		\end{definition}
		Similarly, we can define \textit{fb-continuous, bf-continuous} and \textit{bb-continuous} function. However, if $(Y,\rho)$ is a metric space, then a function $f: (X,d) \to (Y,\rho)$ is ff-continuous (resp. bb-continuous) if and only if it is fb-contiuous (resp. bf-continuous). So in this case, we use the term \textit{forward continuous} for ff-continuous and \textit{backward continuous} for bb-continuous functions.

		For a quasi-metric space $(X,d)$ and nonempty sets $A, B\subseteq X$, we define $d(x,A)= \inf\{d(x,a):a\in A\}$ and $d(A,B)= \inf\{d(a,b):a\in A,\text{ }  b\in B\}$. Recall that if $A$ is a nonempty subset of a metric space $(X,d)$, then the function $f: (X,d) \to (\mathbb{R}, |\cdot|)$ defined by $f(x)=d(x,A)$ is always continuous. However, if $(X,d)$ is a quasi-metric space, then the function $f$ is not necessarily forward continuous.


			We now discuss uniform continuity of a function between quasi-metric spaces.
			\begin{definition}\normalfont Let $(X,d)$ and $(Y,\rho)$ be quasi-metric spaces.
				A function $f:(X,d)\rightarrow(Y,\rho)$ is called \textit{forward uniformly continuous} (resp. \textit{backward uniformly continuous}) if for every $\epsilon >0,$ there exists $\delta>0$ such that $\rho(f(x),f(y))<\epsilon$ (resp. $\rho(f(y),f(x))<\epsilon$) whenever $d(x,y)<\delta$.
			\end{definition}
			It is easy to see that if $(Y,\rho)$ is a metric space, then every forward uniformly continuous function is backward uniformly continuous and vice versa. Consequently, in this case, we simply call such functions uniformly continuous. Need not to mention that like in metric spaces, the notion of uniform continuity depends upon the choice of compatible quasi-metrics on the underlying spaces. 
			   
			\begin{definition}\normalfont (Definition 3.1, \cite{Moshokoa}) Let $(x_n)$ and $(y_n)$ be sequences in a quasi-metric space $(X,d)$. Then $(x_n)$ is said to be \textit{forward parallel} (resp. \textit{backward parallel}) to $(y_n)$ denoted as $ (x_n){\overset{f}{||}} (y_n)$ (resp. $ (x_n)\underset{b}{||} (y_n)$) if for every $ \epsilon >0$, there exists
				$ N \in \mathbb{N}$ such that $d(x_n,y_n) <\epsilon$ (resp. $d(y_n,x_n) <\epsilon$) whenever $ n\geq N$. 
			\end{definition}
		
			The following example shows that these two notions are independent of each other.
			
			\begin{example}
				Let $(\mathbb{R},d_u)$ be the upper quasi-metric space (see, Example \ref{Example upper and lower quasi metrics on R}). Then $x_n=(n+1)$ is forward parallel to $y_n=(n)$ but $(x_n)$ is not backward parallel to $(y_n)$. \qed
			\end{example}
			
			The following proposition provides a characterization of forward uniform continuity, analogous to the characterization in metric spaces, in terms of parallel sequences.
			
			\begin{proposition}\label{p1}$($\cite{Moshokoa}$)$ Let $f:(X,d) \rightarrow (Y,\rho)$ be a function between two quasi-metric spaces $(X,d)$ and $(Y,\rho)$. Then $f$ is forward uniformly continuous if and only if $f$ preserves forward parallel sequences.
			\end{proposition}

			\section{UC Quasi-Metric Spaces}
			In this section, we introduce the notion of a UC quasi-metric space. Several illustrative examples are provided, and the fundamental structural properties of such spaces are examined. It is shown  that, despite some similarities, the behavior of UC quasi-metric spaces exhibits substantial deviations from that of their metric counterparts.  
			
			\begin{definition}\normalfont
				A quasi-metric space $(X,d)$ is called a \textit{UC quasi-metric space} if every forward continuous function $f:(X,d)\rightarrow (\mathbb{R},|\cdot|)$ is uniformly continuous.
			\end{definition} 
			\begin{example} The quasi-metric space $(\mathbb{R}, d_{l})$,  where $d_{l}$ denote the lower quasi-metric as defined in Example \ref{Example upper and lower quasi metrics on R} is a UC quasi-metric space. To see this, we show that every forward continuous function $f : (\mathbb{R}, d_{l}) \to (\mathbb{R}, |\cdot|)$ must be constant. Let $f : (\mathbb{R}, d_{l}) \to (\mathbb{R}, |\cdot|)$ be a forward continuous function.  Suppose there exist two points $x, y \in \mathbb{R}$ such that $f(x) \neq f(y)$. Since the sequence $(x+n)$ converges to both $x$ and $y$ in $(\mathbb{R}, d_{l})$,  it follows that $(f(x+n))$  converges to both $f(x)$ and $f(y)$ in $(\mathbb{R}, |\cdot|)$ . But this is only possible if $f(x) = f(y)$, which contradicts our assumption.	Hence, $f$ must be a constant function. Therefore, $f$ is uniformly continuous. 
				
				In a similar manner, one can see that $(\mathbb{R}, d_{u})$ is UC, where $d_{u}$ denote the upper quasi-metric.\qed
				
			\end{example}
			
		\begin{remark} Note that the associated metric space of $(\mathbb{R}, d_{l})$, that is, $ (\mathbb{R}, |\cdot|)$, is not a UC space.
		\end{remark}
		However, we have the following interesting result. We omit its easy proof. 
		
		\begin{proposition}\label{(X,d) is UC implies (X,d^s) is UC}
			Suppose $(X,d)$ is a $T_1$ quasi-metric space having the property forward convergence implies backward convergence. Then $(X,d^{s})$ is UC whenever $(X,d)$ is UC. 
		\end{proposition}
		As shown by the next example, the converse of the previous result may fail.

		\begin{example}\label{exp of d^s UC but not d uc}
			Let $X = \left\{ \frac{1}{n} : n \in \mathbb{N} \right\} \cup \{0\}$, and let $d$ be the quasi-metric defined as the conjugate of the Sorgenfrey quasi-metric as defined in Example~\ref{Example Sorgenfrey quasi metric}. Then $(X, d)$ is a $T_1$ quasi-metric space and has the property forward convergence implies backward convergence.
			
			Since the associated metric satisfies $d^s(x, y) = 1$ for every $x \neq y$, the space $(X, d^s)$ is a UC metric space. However, $(X, d)$ is not a UC quasi-metric space. To see this, consider the function $f : (X, d) \to (\mathbb{R}, |\cdot|)$, defined by
			
			$$
			f(x) =
			\begin{cases}
				\frac{1}{x}, & \text{if } x \neq 0, \\
				0, & \text{if } x = 0.
			\end{cases}
			$$
			Then $f$ is forward continuous but not uniformly continuous on $(X,d)$.
		\end{example}

		We now show that every compact quasi-metric space is UC. To do this, we first present a more general result, originally due to Lambrinos \cite{Lambrinos1973}, who proved it for quasi-uniform spaces (see also, Theorem 1.1.58 of \cite{Cobzas}). For readers' convenience, we give its full proof in the setting of quasi-metric spaces. We first need the following definition. 
		
		\begin{definition} Let $(X,d)$ be a quasi-metric space, and let $A \subseteq X$. Then $A$ is said to be $\tau(d)$-bounded if every open cover of $X$ in $(X,d)$ contains a finite subcover of $A$.\end{definition} 
		Clearly, $X$ is $\tau(d)$-bounded if and only if $(X,d)$ is compact. Also note that if $A\subseteq X$ is $\tau(d)$-bounded, then for each $\epsilon > 0$ there exists finitely many $x_1, x_2,\ldots,x_n$ in $X$  such that $A\subseteq \cup_{i=1}^{n}B^+(x_i,\epsilon)$, and a closed $\tau(d)$-bounded set is compact. 			
		\begin{theorem} (\cite{Lambrinos1973})
			Let $(X,d)$ be a quasi-metric space. Then every forward continuous function $f: (X,d) \to (\mathbb{R},|\cdot|)$ is uniformly continuous on each $\tau(d)$-bounded subset of $X$.	
		\end{theorem}
		\begin{proof}
			Let $A\subseteq X$ be $\tau(d)$-bounded and $f:(X,d)\rightarrow(\mathbb{R},|\cdot|)$ a forward continuous function. We show that $f:(A,d)\rightarrow (\mathbb{R},|\cdot|)$ is uniformly continuous. Let $\epsilon>0$. As $f:(X,d)\rightarrow (\mathbb{R},|\cdot|)$ is forward continuous, for each $x\in X$, there exists $\delta_x >0$ such that $f(B^{+}(x,\delta_x))\subseteq B(f(x),\frac{\epsilon}{2})$. Then $\{B^{+}(x,\frac{\delta_x}{2}):x\in X\}$ is an open cover of $X$. Since $A$ is $\tau(d)$-bounded, there exists $x_1,x_2,\cdots,x_n \in X$ such that $A \subseteq \cup_{i=1}^{n}B^{+}(x_i,\frac{\delta_{x_{i}}}{2})$. Take $\delta= \min\{\frac{\delta_{x_{i}}}{2}:1\leq i\leq n\}$. Let $x,y\in A$ such that $d(x,y)<\delta$. Then there exists $1\leq k \leq n$ such that $x\in B^{+}(x_k,\frac{\delta_{x_{k}}}{2})$.  So $d(x_k,y)\leq d(x_k,x)+d(x,y)<\delta_{x_{k}}$. Thus, $x,y\in B^{+}(x_k,\delta_{x_{k}})$. Since $f(B^{+}(x_k,\delta_{x_k}))\subseteq B(f(x_k),\frac{\epsilon}{2})$, we get $|f(x)-f(y)|<\epsilon$.	
		\end{proof}
		\begin{corollary}
			If $(X,d)$ is a compact quasi-metric space, then $(X,d)$ is a UC space.
		\end{corollary}
	
		\begin{example}
			The quasi-metric spaces $(\mathbb{R}, d_{l})$ is UC but not compact. As $\{(-\infty, n) : n \in \mathbb{N}\}$ is an open cover of $(\mathbb{R}, d_{l})$ having no finite subcovers.
		\end{example}
		It is well-known that every UC metric space is complete. As mentioned earlier, in the case of quasi-metric spaces, there are several non-equivalent notions of completeness. Among these completeness notions, we relate the left K-completeness with the notion of UC quasi-metric spaces. Our choice is motivated by the fact that left K-complete quasi-metric spaces possess properties analogous to those of complete metric spaces. For instance, every  compact quasi-metric space is left K-complete and every precompact, left K-complete quasi-metric space is compact. Unfortunately, a UC quasi-metric space may not always be left K-complete (see Example \ref{uc does not implies left k-complete}).

		We now examine, when the UC property of a quasi-metric space implies its left K-completeness.

		\begin{theorem}\label{UC implies left k-complete result}
			Suppose $(X,d)$ is a quasi-metric space which is also normal. Then $(X,d)$ is left K-complete, whenever $(X,d)$ is UC.
		\end{theorem}
		\begin{proof} 
			Let $(x_n)$ be any left K-Cauchy sequence in $(X,d)$ which is not convergent. Without loss of generality, we may assume that the terms of the sequence $(x_n)$ are all distinct. Define $A=\{x_{2n}:n\in \mathbb{N}\}$ and $B=\{x_{2n-1}:n\in \mathbb{N}\}$. Since $(x_n)$ cannot not have a convergent subsequence,  $A$ and $B$ are closed subsets of $(X,d)$. Since $(X,d)$ is normal, there exists a forward continuous function  $f:(X,d)\to (\mathbb{R},|\cdot|)$ such that $f(x)=0$ for every $x\in A$ and $f(x)=1$ for every $x\in B$. Therefore, $|f(x_{2n})-f(x_{2n-1})|=1$ for every $n\in \mathbb{N}$. Hence $f$ is not uniformly continuous.
		\end{proof}
		The following example justifies the assumption of normality in Theorem \ref{UC implies left k-complete result}.
		\begin{example}\label{uc does not implies left k-complete}
			
			Let $X=\{\frac{1}{n}:n\geq 2\}\cup \mathbb{N}$ and $d$ be a quasi-metric defined as
			$$
			\begin{aligned}
				&d(x, x) = 0 \quad \text{for all } x \in X,\\[6pt]
				&d\left(\frac{1}{m}, \frac{1}{n}\right) = \frac{1}{m} - \frac{1}{n} \quad \text{if } \frac{1}{m} > \frac{1}{n},\\[6pt]
				&d\left(\frac{1}{m}, n\right) = \frac{1}{n},\\[6pt]
				&d(x, y) = 1 \quad \text{otherwise}.
			\end{aligned}
			$$
			Then $(X,d)$ is a $T_1$ quasi-metric space which is not normal. 
			Note that the sequence $x_n = n$ is forward convergent to $\frac{1}{m}$ for all $m\in \mathbb{N}$ in $(X,d)$. Moreover, $(X,d)$ is not left K-complete as the sequence $(x_n)=(\frac{1}{n})$ is left K-Cauchy in $(X,d)$ but is not forward convergent in $(X,d)$. 
			
			We show that $(X,d)$ is a UC space. Let $f:(X,d)\rightarrow (\mathbb{R},|\cdot|)$ be a forward continuous map. Since for the sequence $x_n = n$ for all $n \in \mathbb{N}$, we have $x_n\xrightarrow{f}\frac{1}{m}$ for all $m\geq 2$, $f(x_n)\rightarrow f(\frac{1}{m})$ for all $m\geq 2$. Consequently, $f(\frac{1}{m})=f(\frac{1}{n})$ for all $m, n \geq 2$. Suppose $f(\frac{1}{n}) = z$ for all $n\geq 2$. To see that $f$ is uniformly continuous, suppose $(x_n)$ is forward parallel to $(y_n)$ in $(X,d)$. We show that $(f(x_n))\parallel (f(y_n))$ in $(\mathbb{R},|\cdot|)$.
			
			First, suppose $(x_n)$ eventually belongs to $\{\frac{1}{n}:n\geq 2\}$, so we can assume $(x_n)$ is a sequence in $\{\frac{1}{n}:n\geq 2\}$. Suppose infinitely many terms of $(y_n)$  are in $\{\frac{1}{n}:n\geq 2\}$ and infinitely many $y_n\in \mathbb{N}$. Then we can divide $(y_n)$ into two subsequences $(y_{m_{k}}), (y_{n_{k}})$ such that $(y_{m_{k}})\subseteq \{\frac{1}{n}:n\geq 2\}$ and $(y_{n_{k}})\subseteq \mathbb{N}$. Since $(x_n)$ is forward parallel to $(y_n)$, for any $\epsilon >0$ there exists $n_0\in \mathbb{N}$ such that  $d(x_{n_{k}},y_{n_{k}})=\frac{1}{y_{n_{k}}}<\epsilon$ for all $ k\geq n_0$. So $d(x_1,y_{n_{k}})=\frac{1}{y_{n_{k}}}<\epsilon$ for all $k\geq n_0$. Hence $y_{n_{k}}\xrightarrow{f}x_1$. Since $f$ is forward continuous, $f(y_{n_{k}})\rightarrow z$. Also $f(y_{m_{k}})=z$ for all $k\in \mathbb{N}$. Therefore, $(f(x_n))\parallel (f(y_n))$. If terms of the sequence $(y_n)$ eventually belongs to $\mathbb{N}$, then a similar argument, as given for $(y_{n_k})$ above, shows that $f(x_n)\parallel f(y_n)$. Finally, if $(y_n)$ eventually belongs to $\{\frac{1}{n}:n\geq 2\}$, then $f(x_n)=f(y_n)=z$ eventually. Hence $(f(x_n))\parallel (f(y_n))$.

			If $(x_n)$ is not eventually from $\{\frac{1}{n}:n\geq 2\}$, then either $(x_n)$ eventually belongs to $\mathbb{N}$ or else we can divide $(x_n)$ into two complementary subsequences $(x_{n_{k}})$ in $\mathbb{N}$ and $(x_{m_{k}})$ in $\{\frac{1}{n}:n\geq 2\}$. In the former case, it is easy to see that $(f(x_n))\parallel (f(y_n))$. In the later case, note that as $(x_n)\parallel (y_n)$, we have $x_{n_{k}}=y_{n_{k}}$ eventually. Moreover, as $(x_{m_{k}})$ is in $\{\frac{1}{n}:n\geq 2\}$, by arguing as above we can show that $(f(x_{m_{k}}))\parallel (f(y_{m_{k}}))$. Consequently, $(f(x_n))\parallel (f(y_n))$. Thus, $(X,d)$ is a UC space. \qed
		\end{example}
		
		Our next result provides a necessary condition for a quasi-metric space $(X,d)$ to be a UC space. Its proof follows from Theorem \ref{Lebesgue Property and UCness} (see, Remark \ref{UC implies X'compact remark}).      
		\begin{theorem}\label{UC implies X' is compact}
			Let $(X,d)$ be a $T_1$ quasi-metric space which is also normal. If $(X,d)$ is a UC space, then the set of non-isolated points $X^{'}$ of $X$ is compact in $(X,d)$.
		\end{theorem}

The following example shows that Theorem \ref{UC implies X' is compact} may not be necessarily true without the assumption of normality.

\begin{example}\label{UC quasi-metric space but X' is not compact}
			Let $(a_n)$ and $(b_n)$ be two sequences of distinct points and such that $a_n\neq b_m$ for any $n, m \in \mathbb{N}$. Let $A = \{ a_n : n \in \mathbb{N} \}$ and $B = \{ b_n : n \in \mathbb{N}\}$. So $A \cap B = \phi$. Let $X = A \cup B$, and define a quasi-metric $d$ on $X$ as follows: 
			$$\begin{aligned}
				&d(a_n, b_m) = \frac{1}{m} \quad \text{for all } n, m \in \mathbb{N},\\	
				&d(a_n, a_m) = 1 \quad \text{ if } n \ne m,\\
				&d(b_n, x) = 1 \quad \text{ if } x \ne b_n, \text{ and} \\
				&d(x, x) = 0 \quad \text{ for all } x \in X.
			\end{aligned}$$
			
			Then $(X,d)$ is a quasi-metric space, and $X'$ is not compact (see, Example 3 of \cite{romagueraAntonino1994}). Note that $(X,d)$ is not normal, since $\{a_{2n}: n \in \mathbb{N}\}$ and $\{a_{2n+1}: n \in \mathbb{N}\}$ are two disjoint closed sets in $(X,d)$ that cannot be separated by open sets. 
			
			We show that $(X,d)$ is a UC space. Let $f:(X,d)\rightarrow (\mathbb{R},|\cdot|)$ be a forward continuous map. Since  $b_n\xrightarrow{f}a_k$ for all $k \in \mathbb{N}$, we have $f(b_n)\rightarrow f(a_k)$ for all $k \in \mathbb{N}$. Consequently, $f(a_n)=f(a_m)$ for all $m,n\in \mathbb{N}$. Let $f(a_n) = z$ for all $n \in \mathbb{N}$. Suppose $(x_n)$ and $(y_n)$ are sequences in $(X,d)$ such that $(x_n)\overset{f}{\parallel}(y_n)$. By a similar argument as given in Example \ref{uc does not implies left k-complete}, we can assume that either $(x_n)$ and $(y_n)$ are eventually equal or $(x_n)$ is eventually in $A$ and $(y_n)$ is eventually in $B$. In the first case, we trivially get $(f(x_n))\parallel (f(y_n))$. In the later case, first choose a $k_0 \in \mathbb{N}$ such that $x_n \in A$ and $y_n \in B$ for all $n \geq k_0$. We claim that $y_n\xrightarrow{f}a_1$. To see this, take an $\epsilon > 0$. Since $(x_n)\overset{f}{\parallel}(y_n)$, we can find an $n_0 \in \mathbb{N}$ such that $n_0\geq k_0$ and $d(x_n,y_n) < \epsilon$ for all $n\geq n_0$. Since $y_n \in B$ for all $n\geq n_0$, let $y_n = b_{m_n}$. Then $\frac{1}{m_n} < \epsilon$ for all $n\geq n_0$. So $d(a_1,y_n) = \frac{1}{m_n} < \epsilon$ for all $n\geq n_0$. Hence $y_n\xrightarrow{f}a_1$. Consequently, $f(y_n) \to f(a_1) = z$. Since each term of the sequence $(f(x_n))$ is eventually equal to $z$, we have $(f(x_n))\parallel (f(y_n))$. Hence $(X,d)$ is a UC quasi-metric space.
		\end{example}

		Our next result is a partial converse of the above theorem. Several proofs are available for a similar result in the context of metric spaces (see, \cite{beer1988b, mrowka1965}). Our proof here is a modification of the proof given by Beer \cite{beer1988b}. A similar result was earlier established in the context of quasi-metric spaces having the Lebesgue number property \cite{Antonioromaguera1993}.

		\begin{theorem}
			Let  $X$ be a quasi-metrizable space and $X^{'}$ be compact. Then there exists a compatible quasi-metric $\rho$ on $X$ such that $(X, \rho)$ becomes a UC space.
		\end{theorem}
		\begin{proof}
			If \( X^{'} = \emptyset \), then the topology of $X$ is discrete. Consequently, the discrete metric $\rho$ on $X$  is a compatible quasi-metric on $X$. Moreover, the space $(X, \rho)$ is a UC space.
			
			So suppose \( X^{'} \neq \emptyset \). Let \( d \) be a compatible quasi-metric for \( X \). Define \( \rho : X \times X \to \mathbb{R} \) as
			\[
			\rho(x, y) = \begin{cases} 
				0 & \text{if } x = y, \\
				d(x, y) + \max \{d(X^{'},x), d(X^{'},y)\} & \text{if } x \neq y.
			\end{cases}
			\]

			One can verify that $\rho$ is quasi-metric on $X$. We show that $d$ and $\rho$ are equivalent. It is easy to see that every forward convergent sequence in $(X, \rho)$ is forward convergent in $(X, d)$. Suppose $(x_n)$ is a  forward convergent sequence  in $(X, d)$ converging to some $x \in X$. If $x \notin X'$, then $(x_n)$ is eventually constant to $x$, and hence $(x_n)$ forward converges to $x$ in $(X,\rho)$. So assume $x \in X^{'}$. Then $d( X^{'},x) = 0$. So $\rho(x,x_n) = d(x,x_n)+d(X^{'},x_n).$ 
			As $d(X^{'},x_n)\leq d( X^{'},x)+ d(x,x_n)$, hence $\rho(x,x_n)\rightarrow 0$. So $(x_n)$ is forward convergent to $x$ in $(X, \rho)$. Therefore, $\rho$ and $d$ are equivalent quasi-metrics on $X$.
			
			Suppose $(X,\rho)$ is not a UC space. Suppose $f:(X,\rho)\rightarrow (\mathbb{R},|\cdot|)$ is a forward continuous function which is not uniformly continuous. So we can find an $\epsilon > 0$ and sequences $(x_n)$ and $(y_n)$ in $(X,\rho)$ such that  $0 < \rho(x_n, y_n) < \frac{1}{n}$, but $|f(x_n)- f(y_n)| > \epsilon$ for all $n\in \mathbb{N}$.  Since $\rho(x_n,y_n)\rightarrow 0$, $d(X^{'},x_n)\rightarrow 0$. Hence we can find a sequence $(t_{k})$ in $X^{'}$ and a subsequence $(x_{n_{k}})$ of $(x_n)$ such that $d (t_{k},x_{n_k})\rightarrow 0$. Since $X^{'}$ is compact, $(t_{k})$ has a cluster point in $(X,d)$, say $a$. Then $a$ is also a cluster point of $(x_{n_{k}})$  in $(X,d)$. Since $d$ and $\rho$ are equivalent, $a$ is a cluster point of $(x_{n_{k}})$ in $(X,\rho)$. So there exists a subsequence $(x_{n_{k_{m}}})$ of $(x_{n_{k}})$ such that $x_{n_{k_{m}}}\xrightarrow{f}a$ in $(X,\rho)$. Consequently, $y_{n_{k_{m}}}\xrightarrow{f}a$. As $f$ is forward continuous, $f(x_{n_{k_{m}}})\rightarrow f(a)$ and $f(y_{n_{k_{m}}})\rightarrow f(a)$, thus we get $|f(x_{n_{k_{m}}})-f(y_{n_{k_{m}}})|<\epsilon$ for all $m\geq n_0$ for some sufficiently large $n_0$. Which contradicts the fact that $|f(x_n)- f(y_n)| > \epsilon$ for all $n\in \mathbb{N}$.
		\end{proof}
		It is to be noted that in the above result, we need not have $\tau(\bar{d}) = \tau(\bar{\rho})$. Moreover, if in the definition of $\rho(x,y)$, we write $d(x,X')$ in place of $d(X',x)$, then $\rho$ and $d$ may not be equivalent quasi-metrics.

		\section{Various Characterizations of UC Quasi-Metric Spaces}				
		The goal of this section is to present several important characterizations of UC quasi-metric spaces. Among others, we characterize these spaces in terms of isolation functionals and the Lebesgue number property. To present the proofs of our main results in a streamlined manner, we first establish an important result related to parallel sequences in quasi-metric spaces.

		Note that in a metric space, if a sequence $(x_n)$ is parallel to a sequence $(y_n)$ and $(x_n)$ does not have a cluster point, then $(y_n)$ also cannot have a cluster point. This fact is helpful in constructing two disjoint closed sets that are near each other. However, this need not be true in a quasi-metric space. In Example \ref{exp of d^s UC but not d uc}, the sequence $(x_n)=(\frac{1}{n})$ is forward parallel to the constant sequence $y_n=0$, and also $(x_n)$ does not have a cluster point. The following result shows that if the underlying quasi-metric space is Hausdorff, one can still construct such closed sets.

		\begin{proposition}\label{existence of disjoint closed sets that are near}
			Let $(X,d)$ be a $T_2$ quasi-metric space. If $(x_n)$ and $(y_n)$ are sequences in $X$ such that $(x_n)$ does not have a cluster point in $(X,d)$ and $ 0<d(x_n,y_n)<\frac{1}{n}$ for all $n \in \mathbb{N}$, then there exist two disjoint closed subsets $C$ and $D$ of $X$ having infinitely many terms of $(x_n)$ and $(y_n)$, respectively such that $d(C,D)=0$.
		\end{proposition}
		
		\begin{proof} Note that the condition $0<d(x_n,y_n)$ for all $n \in \mathbb{N}$ implies that $x_n\neq y_n$ for each $n \in \mathbb{N}$. First, we show that there exist subsequences $(x_{n_{k}})$ and $(y_{n_{k}})$ of $(x_n)$ and $(y_n)$, respectively, such that $x_{n_{k}}\neq y_{n_{l}}$ for every $k,l\in \mathbb{N}$. 
			
			If $(x_n)$ has a constant subsequence $(x_{n_k})$ whose each term is $a$, then take $(y_{n_k})$ to be the corresponding subsequence of $(y_n)$. Similarly, if $(y_n)$ has a constant subsequence $y_{n_k} = b$ for all $k$, then $(x_{n_k})$ be the required subsequence of $(x_n)$.

			So without loss of generality, by passing to subsequences, we may assume that $(x_n)$ and $(y_n)$ are distinct terms sequences. Now proceed as follows:
			
			Take $x_{n_1} = x_1$ and $y_{n_1} = y_1$. Now, choose $x_{n_2}$ and $y_{n_2}$ corresponding parallel elements of $(x_n)$ and $(y_n)$ such that $x_{n_2} \neq y_1$ and $y_{n_2} \neq x_1$.
			Since all the elements of $(x_n)$ and $(y_n)$ are distinct, we can inductively choose $x_{n_{k+1}}$ and $y_{n_{k+1}}$ such that $n_{k+1} > n_k$ and
			$x_{n_{k+1}} \neq y_{n_1}, y_{n_2}, \ldots, y_{n_k} \quad \text{and} \quad y_{n_{k+1}} \neq x_{n_1}, x_{n_2}, \ldots, x_{n_k}$. So $(x_{n_k})$ and $(y_{n_k})$ are the required subsequences.
			
			Now we show the existence of disjoint closed sets. By the above argument, passing to subsequences if necessary, we may assume without loss of generality that $x_n \neq y_m$ for every $n, m \in \mathbb{N}$.
			
			If the sequence $(y_n)$ does not have a cluster point, then take $C = \{x_n : n \in \mathbb{N}\}$ and $D = \{y_n : n \in \mathbb{N}\}$. Suppose $(y_n)$ has a cluster point $a$ in $(X,d)$.
			Then there exists a subsequence $(y_{n_k})$ of $(y_n)$ such that $y_{n_k}\xrightarrow{f} a$. Since the sequence $(y_{n_{k}})$ has a unique limit, $\{y_{n_k} : k \in \mathbb{N}\} \cup \{a\}$ is closed in $(X,d)$. Therefore, $C = \{x_{n_k} : k \in \mathbb{N}\} \setminus \{a\}$ and $D = \{y_{n_k} : k \in \mathbb{N}\} \cup \{a\}$ are the required closed sets. This completes the proof.
		\end{proof}

		We now give a sequential characterization of UC quasi-metric spaces. In order to give this characterization, we first need to define degree of isolation at point $x$ in a quasi-metric space $(X,d)$.  
		
		\begin{definition}\normalfont
			Let $(X,d)$ be a quasi-metric space and $x\in X$. We define $I^{+}(x)=d(x,X\setminus \{x\})= \inf \{d(x,y):y\in X\setminus \{x\}\}$. The functional $I^+$ is called the \textit{isolation functional} on $(X,d)$, and $I^+(x)$ is called the \textit{degree of isolation at point} $x$ in the quasi-metric space $(X,d)$.
			
		\end{definition}
		
		It is easy to see that $I^{+}(x) = 0$ if and only if $x$ is a non-isolated point, that is, $x$ is an accumulation point of $(X,d)$, and that for an isolated point $x$ of $(X,d)$, $I^{+}(x)>0$.

	\begin{theorem}\label{isolation functional characterization}
		Let $(X,d)$ be a $T_1$ quasi-metric space. If $(X,d)$ is normal, then the following assertions are equivalent:
		\begin{enumerate}[(a)]
			\item $(X,d)$ is a UC space;
			\item whenever $(x_n)$ is a sequence in $X$ such that $(I^{+}(x_n))$ converges to $0$, then $(x_n)$ has a cluster point in $(X,d)$;
			\item if $(x_n)$ and $(y_n)$ are two sequences such that $x_n\neq y_n$ for each $n$ and $(x_n)$ is forward parallel to $(y_n)$ in $(X,d)$, then the sequence $(x_n)$ has a cluster point in $(X,d)$.
		\end{enumerate}
	\end{theorem}
	\begin{proof}

		$(a)\implies(b)$. Suppose $(x_n)$ is a sequence in $X$ without a cluster point in $(X,d)$ such that $I^+(x_n) \to 0$. Since $I^{+}(x_n)\rightarrow 0$, we can find a subsequence $(x_{n_k})$ of $(x_n)$ such that $I^{+}(x_{n_{k}})<\frac{1}{k}$ for each $k \in \mathbb{N}$. Consequently, there exists $ y_{n_{k}}\neq x_{n_{k}}$ such that $0<d(x_{n_{k}},y_{n_{k}})<\frac{1}{k}$. As $(x_n)$ has no cluster point, $(x_{n_{k}})$ also has no cluster point in $(X,d)$. By Proposition \ref{existence of disjoint closed sets that are near}, we get two disjoint closed subsets $C$ and $D$ containing infinitely many points of $(x_{n_{k}})$ and $(y_{n_{k}})$, respectively such that $d(C,D)=0$. Since $(X,d)$ is normal, there exists a forward continuous function $f : (X, d) \to (\mathbb{R},|\cdot|)$ such that $f(C) = 0$ and $f(D) = 1$. Then  $(f(x_{n_{k}})){\not\parallel} (f(y_{n_{k}}))$, while $(x_{n_{k}}) \overset{f}{\parallel} (y_{n_{k}})$. It follows that $f$ is not uniformly continuous.
		
		$(b)\implies (c)$. Suppose $(x_n)$  is forward parallel to $(y_n)$ in $(X,d)$ and $x_n\neq y_n$ for every $n\in \mathbb{N}$. Then it is easy to see that $(I^{+}(x_n))$ converges to $0$, and thus $(x_n)$ has a cluster point.

		$(c) \implies (a)$. Suppose $(X,d)$ is not a UC space. Then there exists a forward continuous function $f:(X,d)\to (\mathbb{R},|\cdot|)$ which is not uniformly continuous. So there exists an $\epsilon > 0$ such that for each $n \in \mathbb{N}$, there exist points $x_n$ and $y_n$ in $X$ with $0 < d(x_n, y_n) < \frac{1}{n}$, but $|f(x_n)- f(y_n)| > \epsilon$. Therefore, $(x_n)\overset{f}{||}(y_n$). By the hypothesis in $(c)$, the sequence $(x_n)$ must cluster at some point $x$ in $(X,d)$. Let $(x_{n_k})$ be a subsequence of $(x_n)$ such that $x_{n_k} \xrightarrow{f} x$. Then $y_{n_k} \xrightarrow{f} x$ in $(X,d)$. Since $f$ is forward continuous, $f(x_{n_k}) \rightarrow f(x)$ and $f(y_{n_k}) \rightarrow f(x)$ in $(\mathbb{R},|\cdot|)$. Therefore, for sufficiently large $k$, we must have
		\[
		|f(x_{n_k})- f(y_{n_k})| < |f(x_{n_k})- f(x)| + |f(x)- f(y_{n_k})|<\epsilon.
		\]
		But this contradicts the fact that $|f(x_{n_k})-f(y_{n_k})| > \epsilon$ for all $k$. Hence $(X,d)$ must be a UC space.
	\end{proof}
	\begin{remark} \label{chr of uc space in terms of degree of isolation functional remark}
		We are using normality of $(X,d)$ only in the proof of the implication $(a)\Rightarrow (b)$.
	\end{remark} 
	
	The following example illustrates why it is necessary for the space $(X, d)$ to be normal.
	\begin{example} \label{counterexample of UC implies isolation functional characterization}
		Let $X=(0,\infty)$ and 	$d :(0,\infty) \times (0,\infty) \rightarrow  \mathbb{R}$ defined by
		$$ d(x,y)=\begin{cases}
			0, \text{ if } x=y\\    
			y, \text{   if }  y\in \{\frac{1}{n}:n\in \mathbb{N}\}\\
			1, \text{ otherwise.}
		\end{cases}$$ 
		
		Then $(X, d)$ is a $T_1$ quasi-metric space. For any fixed $m \in \mathbb{N}$, the sequence $\left(\frac{1}{n}\right)$ forward convergent to $\frac{1}{m}$. Consequently, $(X,d)$ is not $T_2$, and hence not normal.  Since every forward continuous map $f : (X, d) \to (\mathbb{R}, |\cdot|)$ is constant, $(X, d)$ is a UC space. However, the sequence $x_n = n$ does not cluster in $(X, d)$, while $I^{+}(x_n) \to 0$. \qed
	\end{example}
	In order to prove the next result, first we need the following definition.
	\begin{definition}\normalfont
		Let $(X, d)$ be a quasi-metric space and $\phi\neq A \subseteq X$. Then $A$ is called \textit{uniformly discrete} if there exists a $\delta > 0$ such that
		\[
		d(x, y) \geq \delta \quad \text{for all} \quad x, y \in A, x\neq y.
		\]
	
	\end{definition}
	Note that $A$ is uniformly discrete in $(X,d)$ actually implies $A$ is uniformly discrete in $(X,d^s)$. However, the converse may not be true. Also a uniformly discrete set in $(X,d)$ is always discrete.
	
	In the next theorem, we are assuming a condition stronger than normality on $(X,d)$. However, we would like to mention that $(a)\Leftrightarrow (b)$ and $(b)\Rightarrow (c)$ in the next theorem can also be proved when $(X,d)$ is normal.  
	\begin{theorem} \label{UCness and uniform discreteness}
		Let $(X,d)$ be a $T_1$ quasi-metric space having the property forward convergence implies backward convergence. Then the following statements are equivalent:
		\begin{enumerate}[(a)]
			\item $(X,d)$ is a UC space;
			\item if $A$ and $B$ are two nonempty disjoint closed subsets of $X$, then $d(A,B)>0$;
			\item  every closed discrete subset of $(X,d)$ is uniformly discrete.
		\end{enumerate}
	\end{theorem}
	\begin{proof}
		$(a)\implies (b)$.	Suppose on the contrary that there exists two nonempty disjoint closed subsets $A_1$ and $A_2$ of $X$ such that $d(A_1,A_2)=0$. So there exists a forward continuous map $f:(X,d)\rightarrow (\mathbb{R},|\cdot|)$ such that $f(A_1)=0$ and $f(A_2)=1$. Since $d(A_1,A_2)=0$, we can find sequences $(x_n)$ and $(y_n)$ in $A_1$ and $A_2$, respectively such that $d(x_n,y_n) < \frac{1}{n}$ for each $n$. However, $(f(x_n))\overset{f}{\not \parallel}(f(y_n))$. We arrive at a contradiction.
		
		$(b)\implies (c)$.	Suppose there exists a closed discrete subset $Y$ of $(X,d)$ which is not uniformly discrete. Then for every $n \in \mathbb{N}$, there exist distinct points $x_n, y_n \in Y$ such that $0 < d(x_n, y_n) < \frac{1}{n}$. If $(x_{n_k})$ is a constant subsequence of $(x_n)$, say $x_{n_k} = a$ for all $k$, then $(y_{n_k})$ forward converges to $a$. But then $(y_{n_k})$ must be eventually constant as $Y$ is closed and discrete. Which contradicts $0 < d(x_n, y_n) < \frac{1}{n}$ for all $n \in \mathbb{N}$. So we can assume $(x_n)$ has distinct terms. Therefore, $(x_n)$ cannot have a cluster point. Then by Proposition \ref{existence of disjoint closed sets that are near}, we can find disjoint closed sets $C$ and $D$ in $(X,d)$ such that $d(C,D)  =0$. We arrive at a contradiction.

		$(c)\implies (a)$. Suppose $(X,d)$ is not a UC space. Then there exists a forward continuous function $f: (X,d) \to (\mathbb{R},|\cdot|)$ which is not uniformly continuous. Therefore, we can find an $\epsilon > 0$ and sequences $(x_n)$ and $(y_n)$ in $(X,d)$ such that $x_n\neq y_m$ for every $n,m \in \mathbb{N}$ and $0 < d(x_n, y_n) < \frac{1}{n}$, but $|f(x_n)- f(y_n)| > \epsilon$. First we show that $(x_n)$ does not have a cluster point in $(X,d)$. If $(x_n)$ has a cluster point $x$ in $(X,d)$, then there exists a subsequence $(x_{n_{k}})$ of $(x_n)$ such that $x_{n_{k}}\xrightarrow{f} x$. Consequently, $y_{n_{k}}\xrightarrow{f} x$. Since $f$ is forward continuous, $f(x_{n_{k}})\rightarrow f(x)$ and $f(y_{n_{k}})\rightarrow f(x)$. Therefore, there exists $k_0\in \mathbb{N}$ such that $| f(x_{n_{k}})-f(y_{n_{k}})|<\epsilon$  $ \forall k\geq k_0$, a contradiction.      
		
		Now we show that the sequence $(y_n)$ has no cluster point in $(X,d)$.
	
		Suppose $(y_n)$ clusters to some point $a$ in $(X,d)$. Take $C = \{x_n : n \in \mathbb{N}\} \cup \{a\}$.  Since $a$ is not a cluster point of $(x_n)$, it is an isolated point of $C$, and hence $C$ is closed and discrete. Thus, $C$ is uniformly discrete. So there is a $\delta>0$ such that $d(x,y)\geq \delta$ for all $x, y \in C$ with $x\neq y$. Since $(y_n)$ clusters to $a$, there exists a subsequence $(y_{n_k})$ of $(y_n)$ such that $y_{n_k} \underset{b}{\overset{f}{\rightarrow}} a$. Therefore, we can find $k_0\in\mathbb{N}$ satisfying $\frac{1}{k_0}<\frac{\delta}{2}$, $d(a, y_{n_{k_0}}) < \frac{\delta}{2}$, $d(y_{n_{k_0}}, a) < \frac{\delta}{2}$, and $x_{n_{k_0}}\neq a$.  It follows that
		\[
		d(x_{n_{k_0}}, a) \leq d(x_{n_{k_0}}, y_{n_{k_0}}) + d(y_{n_{k_0}}, a) < \frac{\delta}{2} + \frac{\delta}{2} = \delta.
		\]
		This contradicts the choice of $\delta >0$.  Therefore, $(y_n)$ does not have a cluster point in $(X,d)$.
		
		Then $A = \{x_n : n \in \mathbb{N}\} \cup \{y_n : n \in \mathbb{N}\}$ is closed and discrete, but not uniformly discrete in $(X,d)$, since $0 < d(x_n, y_n) < \frac{1}{n}$ for every $n \in \mathbb{N}$.	
	\end{proof}
	
	Example \ref{uc does not implies left k-complete} shows that in the absence of normality the equivalence $(a)\Leftrightarrow (b)$ in the above result may fail. For instance in Example \ref{uc does not implies left k-complete}, if we take $A = \{\frac{1}{2n}: n \in \mathbb{N}\}$ and $B = \{\frac{1}{2n+1}: n \in \mathbb{N}\}$, then $A$ and $B$ are disjoint closed sets in $(X,d)$ such that $d(A,B) =0$.

	\begin{remark}
		The quasi-metric spaces satisfying the condition $(b)$ of Theorem \ref{UCness and uniform discreteness} are sometimes called equi-normal in the literature (see, page 70 of \cite{Cobzas}). 
	\end{remark} 
	
	\begin{definition}\label{Lebesgue number property defi}\normalfont
		Let $(X, d)$ be a quasi-metric space, and let $\mathcal{U}$ be an open cover of $X$. A $\delta > 0$ is called a \textit{Lebesgue number} for $\mathcal{U}$ if for every point $x \in X$, the open ball $B^{+}(x, \delta) = \{ y \in X : d(x, y) < \delta \}$ is contained in some $U \in \mathcal{U}$. A quasi-metric space $(X,d)$ is said to have the Lebesgue number property provided every open cover of $X$ has a Lebesgue number.
	\end{definition}
	The Proposition 2.1 in \cite{MarinRomaguera1996} implies that if a quasi-metric space $(X,d)$ has the Lebesgue number property, then every forward continuous function from $(X,d)$ to any quasi-metric space $(M,\rho)$ with the property forward convergence implies backward convergence is forward uniformly continuous. In particlular, it follows that every quasi-metric space with the Lebesgue number property is a UC space.  In the next result, we give a different proof for this fact using isolation functional, and also see when the converse is true.
	\begin{theorem} \label{Lebesgue Property and UCness}
		Let $(X,d)$ be a $T_1$ quasi-metric space that is also normal. Then the following statements are equivalent:
		\begin{enumerate}[(a)]
			\item $(X,d)$ is a UC space;
			\item every open cover of $(X,d)$ has a Lebesgue number.
		\end{enumerate}
	\end{theorem}
	\begin{proof}		
		$(a)\implies(b)$. Let $\{O_{\lambda}:\lambda\in \Lambda\}$ be an open cover of $X$ having no Lebesgue number. So for every $n \in \mathbb{N}$, there exists $x_n \in X$ such that $B^{+}(x_n, \frac{1}{n}) \not\subseteq O_{\lambda}$ for any $\lambda \in \Lambda$. Consequently, for each $n\in \mathbb{N}$, there exists $y_n\in B^{+}(x_n, \frac{1}{n})$ such that $x_n\neq y_n$. Therefore $0<d(x_n,y_n)< \frac{1}{n}$ $\forall n\in \mathbb{N}$.  
		We now show that $(x_n)$ does not have a cluster point in $(X,d)$. If $(x_n)$ has a cluster point, say $x$, in $(X,d)$, then there exists a subsequence $(x_{n_k})$ forward converging to $x$. Since $\{O_{\lambda}:\lambda\in \Lambda\}$ is an open cover of $X$, there exists $\lambda \in \Lambda$ such that $B^{+}(x, \epsilon) \subseteq O_{\lambda}$ for some $\epsilon > 0$. Since $(x_{n_k})$ forward converging to $x$, for sufficiently large $k$, we have $B^{+}(x_{n_k}, \frac{1}{n_k}) \subseteq B^{+}(x, \epsilon)\subseteq O_{\lambda}$. Which contradicts the fact that $B^{+}(x_n, \frac{1}{n}) \not\subseteq O_{\lambda}$ for every $\lambda \in \Lambda$. So $(x_n)$ does not have a cluster point in $(X,d)$. Therefore, by Proposition \ref{existence of disjoint closed sets that are near}, we get two closed and disjoint sets $C$ and $D$ containing infinitely many points of $(x_n)$ and $(y_n)$, respectively, such that $d(C,D)=0$. Which contradicts Theorem \ref{UCness and uniform discreteness}.
		
		$(b)\implies (a)$.
		Suppose $(X,d)$ is not a UC space. Then by Theorem \ref{isolation functional characterization}, we can find a sequence $(x_n)$ in $X$ such that $I^{+}(x_n) \to 0$, while $(x_n)$ does not have a cluster point in $(X,d)$. Without loss of generality, we can assume $(x_n)$ has distinct terms. Since $I^{+}(x_n) \to 0$, we can find a subsequence $(x_{n_k})$ of $(x_n)$ and a sequence $(y_k)$ in $X$ such that $y_k \neq x_{n_k}$ for all $k$ and $0 < d(x_{n_k}, y_k) = \epsilon_k < \frac{1}{k}$. Then $(x_{n_k})$ also cannot have a cluster point in $(X,d)$. Consequently, the set $A = \{x_{n_k} : k \in \mathbb{N}\}$ is a closed and discrete subset of $(X,d)$. Therefore, for every $x_{n_m} \in A$, there exists $\delta_m > 0$ such that $B^{+}(x_{n_m}, \delta_m) \cap A = \{x_{n_m}\}$.
		
		Let $\gamma_m = \min\{\epsilon_m, \delta_m\}$ for each $m\in\mathbb{N}$. Consider the collection 
		$$
		\{B^{+}(x_{n_m}, \gamma_m):m\in \mathbb{N}\}\cup \{X\setminus A\},
		$$
		which is an open cover of $X$. We claim that this cover does not have a Lebesgue number. Suppose on the contrary that $r>0$ is a Lebesgue number for this cover. Choose $k \in \mathbb{N}$ such that $\epsilon_k < \frac{1}{k} < r$. Hence $\gamma_k < \frac{1}{k} < r$. We show that $B^{+}(x_{n_k}, r)$ is not contained in any member of the cover.
		
		If $B^{+}(x_{n_k}, r) \subseteq B^{+}(x_{n_m}, \gamma_m)$ for some $m \neq k$, then $d(x_{n_m}, x_{n_k}) < \gamma_m < \delta_m$. Which contradicts the choice of $\delta_m$. Also, $B^{+}(x_{n_k}, r) \nsubseteq B^{+}(x_{n_k}, \gamma_k)$ as $y_k \in B^{+}(x_{n_k}, r)$ but $y_k \notin B^{+}(x_{n_k}, \gamma_k)$.  Thus, $B^{+}(x_k, r)$ is not contained in any member of the cover $\{B^{+}(x_{n_m}, \gamma_m):m\in \mathbb{N}\}\cup \{X\setminus A\}$. 
	\end{proof}
		
	We now give an example of a non-normal, $T_1$ quasi-metric space $(X,d)$ which is UC but does not satisfy the Lebesgue number property. 
	
	\begin{example}
		Let $(X, d)$ be the quasi-metric space defined in Example \ref{counterexample of UC implies isolation functional characterization}. Consider the open cover 
		\[
		\mathcal{U} = \left\{ B^{+}\left(n, \frac{1}{n}\right) : n \in \mathbb{N} \right\} \cup \left\{ B^{+}\left(x, \frac{1}{2}\right) : x \neq n \text{ and } x \neq \frac{1}{n} \text{ for any } n \in \mathbb{N} \right\}.
		\]
		Suppose, by contradiction, that there exists $\delta > 0$ which serves as a Lebesgue number for this cover. Then there exists $n_0 \in \mathbb{N}$ such that $\frac{1}{n_0} < \delta$. However, $B^{+}(n_0, \delta) \nsubseteq B^{+}(n_0, \frac{1}{n_0})$, and it also cannot be contained in any other member of the cover. Indeed, $n_0 \in B^{+}(n_0, \delta)$, but $n_0 \notin B^{+}(x, \frac{1}{2})$ for any $x \in X$ with $x \neq n_0$. Hence the open cover $\mathcal{U}$ does not admit a Lebesgue number, even though the space is a UC space.\qed
	\end{example}
	
\begin{remark}\label{UC implies X'compact remark}
In \cite{romagueraAntonino1990,romagueraAntonino1994}), Romaguera and Antonino proved that if a quasi-metric space $(X,d)$ has the Lebesgue number property, then the set $X'$ of non-isolated points of $(X,d)$ is compact. Consequently, by Theorem \ref{Lebesgue Property and UCness}, one can deduce that if a quasi metric space $(X,d)$ is $T_1$, normal and UC, then $X'$ is compact (see, Theorem \ref{UC implies X' is compact}). 	
\end{remark}
	Recall that for a metric space $(X,d)$, the following statements are equivalent:
	\begin{enumerate}
		\item every continuous function $f: (X,d) \to (\mathbb{R},|\cdot|)$ is uniformly continuous;
		\item for any metric space $(Y,\rho)$, every continuous function $f: (X,d) \to (Y,\rho)$ is uniformly continuous.
	\end{enumerate}      
	However, we do not have the same privilege in the context of a quasi-metric space. Our next result prove a similar result under certain restrictions on the underlying quasi-metric space.
	\begin{theorem}\label{UC iff strong UC}
		Let $(X,d)$ be a $T_1$ quasi-metric space. If $(X,d)$ is normal, then the following assertions are equivalent:
		\begin{enumerate}[(a)]
			\item $(X,d)$ is a UC space;
			\item if $A$ is a subset of $X$ without an accumulation point in $(X,d)$, then $A\cap X^{'}$ is finite and  $\inf\{I^{+}(x):x\in A\setminus X^{'}\}>0$;
			\item whenever $(x_n)$ is a sequence in $X$ without a cluster point and $B=\{n\in \mathbb{N}:x_n$ is not an isolated point in  $(X,d)\}$, then $B$ is finite and $\inf\{I^{+}(x_n):n\notin B\}>0$; 
			\item for any quasi-metric space $(M,\rho)$ with the property forward convergence implies backward convergence, every forward continuous function $f:(X,d)\rightarrow (M,\rho)$ is uniformly continuous.
		\end{enumerate}
	\end{theorem}
	
	\begin{proof} $(a)\implies (b)$. Let $A$ be a subset of $X$ without an accumulation point in $(X,d)$. First, we show that $A\cap X'$ is finite. By contradiction, suppose $A\cap X^{'}$ is infinite. So there exists a sequence $(x_n)$ of distinct points in $A\cap X^{'}$. Then $(x_n)$ cannot have a cluster point in $(X,d)$ and for each $n\in \mathbb{N}$, there exists $y_n\in X$ such that $0<d(x_n,y_n)<\frac{1}{n}$.
		 
		By Proposition \ref{existence of disjoint closed sets that are near}, there exists two disjoint closed subsets $C$ and $D$ of $X$ having infinitely many terms of $(x_n)$ and $(y_n)$, respectively such that $d(C,D)=0$. Since $(a)\Leftrightarrow (b)$ of Theorem \ref{UCness and uniform discreteness} is true for normal $T_1$ space, we get $(X,d)$ is not a UC space. Hence $A\cap X^{'}$ must be finite. 
		
		Now we show $\inf\{I^{+}(x_n):x\in A\setminus X^{'}\}>0$.  If $A\setminus X^{'}$ is finite, then clearly  $\inf\{I^{+}(x):x\in A\setminus X^{'}\}=\min\{I^{+}(x):x\in A\setminus X^{'}\}>0$. So assume $A\setminus X^{'}$ is infinite. Suppose $\inf\{I^{+}(x):x\in A\setminus X^{'}\}=0$. So we can find a sequence $(x_n)$ of distinct terms in $A\setminus X'$ and a sequence $(y_n)$ in $X$ such that $0<d(x_n,y_n)<1/n$. Then $I_{+}(x_n)\rightarrow 0$ but $(x_n)$ does not have a cluster point in $(X,d)$. Consequently, by Theorem \ref{isolation functional characterization}, we get $(X,d)$ is not UC. Hence $\inf\{I^{+}(x):x\in A\setminus X^{'}\}>0$.
		
		$(b)\implies (c)$. It is immediate.

		$(c)\implies (d)$. 	 Let $(M,\rho)$ be a quasi-metric space with the property forward convergence implies backward convergence. Suppose $f:(X,d)\rightarrow(M,\rho)$ is a forward continuous function which is not uniformly continuous. So there exist sequences $(x_n)$, $(y_n)$ in $X$ and $\epsilon>0$ satisfying $0<d(x_n,y_n)<\frac{1}{n}$ and $\rho(f(x_{n}),f(y_{n}))>\epsilon$ for every $n\in\mathbb{N}$. Then $(x_n)$ does not have a cluster point in $(X,d)$.  Then $B=\{n\in \mathbb{N}:x_n \text{ is not an isolated } \text{ point in } (X,d)\}$ is finite but $\inf\{I^{+}(x_n):n\notin B\}=0$.

		$(d)\implies (a)$. It is immediate.
	\end{proof}

If we drop the condition forward convergence implies backward convergence in the previous result, then it may not be true. The following Example illustrates this.
	
	\begin{example}
		Let $X=\{n:n\in \mathbb{N}\}\cup \{0\}$ along with the quasi-metric $d$ defined as
		\[
		d(x, y) = 
		\begin{cases} 
			d(x, 0) + d(0, y), & \text{if } y \neq x \\
			0, & \text{if } x = y,
		\end{cases}
		\]
		where $d(x, 0) = \frac{1}{x}$ and $d(0, x) = \frac{1}{x^2}$ for $x \neq 0$. Then $(X,d)$ is a metrizable quasi-metric space (see, \cite{Zimmer}). Since for any sequence $(x_n)$ in $X$ such that $I_{+}(x_n)\rightarrow 0$,  $(x_n)$ must be either an eventually constant sequence or has an increasing subsequence. Therefore, whenever $I_{+}(x_n)\rightarrow 0$, $(x_n)$ clusters in $(X,d)$. So $(X,d)$ is a UC space. Let $M=\{\frac{1}{n}\}\cup \{0\}$, equipped with the Sorgenfrey quasi-metric $\rho$. Then $(M,\rho)$ is a $T_1$ quasi-metric but it does not have the property forward convergence implies backward convergence. Define $f:(X,d)\rightarrow (M,\rho)$ as  $$
		f(x) =
		\begin{cases}
			\frac{1}{x}, & \text{if } x \neq 0, \\
			0, & \text{if } x = 0.
		\end{cases}
		$$
		Then $f$ is forward continuous but not uniformly continuous, since $x_n=n$ is forward parallel to $y_n=n+1$ but $(f(x_n))$ is not forward parallel to $(f(y_n))$. \qed
	\end{example}
	\section{More on UC Quasi-Metric Spaces}		
	In this section, we find a Cantor intersection type characterizations for UC spaces. In order to get this, first we need the following definition and result.
	\begin{definition}\normalfont
		A sequence $(x_n)$ in a quasi-metric space $(X,d)$ is called \textit{pseudo left K-Cauchy} (resp. \textit{pseudo right K-Cauchy}) if for every $\epsilon>0$ and $n \in \mathbb{N}$, there exist $j,k\in\mathbb{N}$ such that $d(x_j,x_k)<\epsilon$ (resp. $d(x_k,x_j)<\epsilon$) for $k>j>n$.

	\end{definition}

	\begin{theorem}\label{pseudo cauchy}
		Let $(X,d)$ be a $T_1$ and normal quasi-metric space such that $\tau(d)\subseteq \tau(\bar{d})$. Then the following statements are equivalent:
		\begin{enumerate}[(a)]
			\item $(X,d)$ is a UC space;
			\item every pseudo left K-Cauchy sequence in $(X,d)$ with distinct terms has a cluster point.				
		\end{enumerate}
	\end{theorem}
	\begin{proof}
		$(a) \implies (b)$.  
		Let $(x_n)$ be a pseudo left $K$-Cauchy sequence with distinct terms.  
		So we can find subsequences $(x_{n_k})$ and $(y_{m_k})$ of $(x_n)$ such that $0<d(x_{n_k}, y_{m_k}) < \frac{1}{k}$. Hence $I^{+}(x_{n_k}) \rightarrow 0$. Since $(X,d)$ is a UC space, by Theorem \ref{isolation functional characterization}, the sequence $(x_{n_k})$ clusters. Therefore, $(x_n)$ also clusters.
		
		$(b) \implies (a)$.  
		Suppose that $(X,d)$ is not a UC space.  
		By Theorem \ref{UCness and uniform discreteness}, there exists a closed discrete subset $A \subseteq X$ which is not uniformly discrete.  
		Thus, for every $n \in \mathbb{N}$, there exist $x_n,y_n\in A$ such that $0 < d(x_n, y_n) < \frac{1}{n}$. Since $A$ is closed and discrete in $(X,d)$, by following the similar steps as in the proof of $(b) \implies (c)$ in Theorem \ref{UCness and uniform discreteness}, one can assume that $(x_n)$ has distinct terms. Since $\tau(d)\subseteq\tau(\bar{d})$,  $A$ is closed and discrete in $(X,\overline{d})$. So we can also assume  $(y_n)$ has distinct terms. Hence $(x_n)$ and $(y_n)$ cannot have a cluster point in $(X,d)$. By Proposition \ref{existence of disjoint closed sets that are near}, we can find subsequences $(x_{n_k})$ and $(y_{n_k})$ of $(x_n)$ and $(y_n)$, respectively such that $x_{n_k}\neq y_{n_l}$ for every $k,l \in\mathbb{N}$ and $d(x_{n_k},y_{n_k})<\frac{1}{k}$. Then the sequence $(z_k) = (x_{n_1}, y_{n_1}, x_{n_2}, y_{n_2}, \ldots)$ is a pseudo left $K$-Cauchy sequence in $(X,d)$ with distinct terms which has no cluster point in $(X,d)$. 
	\end{proof}

	The following example demonstrates that if we drop the condition $\tau(d)\subseteq\tau(\bar{d})$ in the previous theorem, then the implication $(b) \implies (a)$ may fail.
	\begin{example}\label{uc implies every pseudo left k cauchy clusters}
		Let $X=\{\frac{1}{n}:n\in \mathbb{N}\}\cup \{0\}$ and $d$ be a quasi-metric defined as 
		
		$$
		\begin{aligned}
			&d(x, x) = 0 \quad \text{for all } x \in X,\\[6pt]
			&d\left(0, \frac{1}{n}\right) = 1,\\[6pt]
			&d\left(\frac{1}{n}, 0\right) = \frac{1}{n},\\[6pt]
			&d(x, y) = 1 \quad \text{otherwise}.
		\end{aligned}
		$$
		Then $(X,d)$  is $T_1$ quasi-metric space having the property forward convergence implies backward convergence. However, $(X,d)$ does not have the property backward convergence implies forward convergence. Moreover, $(X,d)$ does not have any pseudo left K-Cauchy sequence of distinct terms. However, $(X,d)$ is not a UC space, since $f:(X,d)\rightarrow (\mathbb{R},|\cdot|)$ defined by $$
		f(x) =
		\begin{cases}
			\frac{1}{x}, & \text{if } x \neq 0, \\
			0, & \text{if } x = 0
		\end{cases}
		$$ is forward continuous but not uniformly continuous, as $(x_n)=\frac{1}{n}$ is forward parallel to $(y_n=0)$ but $(f(x_n))\nparallel (f(y_n))$. \qed	
	\end{example}
	The following example demonstrates that in the absence of normality of $(X,d)$, the implication $(a) \implies (b)$ may fail in the previous theorem. 
	\begin{example}
		Let $(X,d)$ be the quasi-metric space as defined in Example \ref{uc does not implies left k-complete}. Then $(X,d)$  is a non-normal, $T_1$, UC quasi-metric space having the property backward convergence implies forward convergence. We show that the sequence $(x_n)=(\frac{1}{n+1})$ is a pseudo left K-Cauchy sequence in $(X,d)$ which does not cluster.  If $(x_n)$ clusters in $(X,d)$, then it cannot cluster at any $k\in \mathbb{N}$ as $d(k,\frac{1}{n})=1$ for all $n\in \mathbb{N}$. Moreover, for any $k\geq 2$,  we have $d(\frac{1}{k},x_{n})\geq\frac{1}{k}-\frac{1}{k+1}$ for all $n\in \mathbb{N}$. Hence $(x_n)$ cannot cluster at any point of $\{\frac{1}{n}:n\geq 2\}$.\qed
	\end{example}
	
Let $(X,d)$ be a quasi-metric space, and $\emptyset \neq F \subseteq X$. For any sequence $(s_n)$ in $F$, define  
		$$V_F((s_n))=\sup\{d(s_{2n-1},s_{2n}):n\in \mathbb{N}\}.$$
		
\begin{theorem} [Cantor-Intersection Type Characterization] Let $(X,d)$ be a $T_1$ and normal quasi-metric space such that $\tau(d)\subseteq\tau(\bar{d})$. Then the following statements are equivalent:			\begin{enumerate}
\item  $(X,d)$ is a UC space;
\item for every decresing sequence $F_1\supseteq F_2\supseteq \dots$ of nonempty closed subsets of $X$, if for every $n\in\mathbb{N}$, there exists a sequence $(s_{m}^{n})\in F_n$ such that $V_{F_{n}}((s_{m}^{n}))\neq 0$  and $ \lim_{n\rightarrow \infty}V_{F_n}((s_{m}^{n}))=0$, then $ \bigcap_{n=1}^{\infty} F_n \neq \phi$.
\end{enumerate}
\end{theorem}
\begin{proof}
		
$(a)\implies (b)$ Let $(F_n)$ be any decreasing sequence of nonempty closed subsets of $X$ and $(s_{m}^{n})\in F_n$ such that $V_{F_{n}}((s_{m}^{n}))\neq 0 $ $\forall n$ and $\lim_{n\rightarrow \infty}V_{F_n}((s_{m}^{n}))=0$. Take $x_n = s_{2m-1}^{n}$ such that $s_{2m-1}^{n} \neq s_{2m}^{n}$ for $m \in \mathbb{N}$. Then 
$I^{+}(x_n) = d(x_n, X \setminus \{x_n\}) \leq d(s_{2m-1}^{n}, s_{2m}^{n}) \leq V_{F_{n}}((s_{m}^{n}))$.  
As $\lim_{n \rightarrow \infty} V_{F_n}((s_{m}^{n})) = 0$, it follows that $I^{+}(x_n) \rightarrow 0$. By Theorem~\ref{isolation functional characterization}, $(x_n)$ has a cluster point $x$. Then there exists a subsequence $(x_{n_k})$ of $(x_n)$ such that $x_{n_{k}} \xrightarrow{f} x$. Since for each $n\in \mathbb{N}$, $x_n\in F_{n}$ and $(F_{n})$ is decreasing, $x_{n_k} \in F_m$ for every $k \geq m$. So by closedness of $F_m$, we get $x \in F_m$ for every $m \in \mathbb{N}$. Hence $x \in \bigcap_{n=1}^{\infty} F_n.$

$(b)\implies (a)$ Let $(x_n)$ be a pseudo left K-Cauchy sequence of distinct terms  in $X$. So we can find a subsequence $(z_n)$ of $(x_n)$ such that $d(z_{2n-1},z_{2n})<\frac{1}{n}$.  For each $n$, let $A_n= \overline{\{z_k: k\geq 2n-1\}}$, and $(s_{m}^{n})=(z_{2n-1},z_{2n},z_{2n+1},\ldots)$. So $V_{A_n}((s_{m}^{n})) < \frac{1}{n}$. Consequently, $\lim _{n\rightarrow \infty} V_{A_n}((s_{m}^{n}))= 0$. As $(z_n)$ is a sequence of distinct terms, $V_{A_n}((s_{m}^{n}))\neq 0$ for each $n\in \mathbb{N}$ and $(A_n)$ is a decreasing sequence of closed sets. Hence $\bigcap_{n=1}^{\infty} A_n \neq \emptyset$. Let $x\in \bigcap_{n=1}^{\infty} A_n$. We show that $x$ is a cluster point of $(z_n)$, consequently of $(x_n)$. Choose any $\epsilon >0$ and $n\in \mathbb{N}$. Since $x\in A_n$, there exists $k\in \mathbb{N}$ such that $k>n$ and  $d(x,z_{n_{k}})<\epsilon$. Therefore, $x$ is a cluster point of the sequence $(z_n)$ in $(X,d)$. By Theorem \ref{pseudo cauchy}, $(X,d)$ is a UC space.
	\end{proof}
	
		
We end this paper, by presenting a characterization of UC quasi-metric spaces, given in terms of non-vanishing uniformly continuous real-valued functions. Here we would like to mention that in \cite{beerGarridoMerono2018} Beer and Garrido proved a similar result in the realm of metric spaces. 
		
We leave the proof of the following result as an easy exercise. 
		
\begin{proposition} \label{closed subspace of a UC space is UC}
If $(X,d)$ is a $T_1$ quasi-metric space which is also normal, then every closed subset of a UC-space is itself a UC-space.
\end{proposition}
		
\begin{theorem}
Let $(X,d)$ be a quasi-metric space such that $\tau(d)=\tau(\bar{d})$. Then $(X,d)$ is a UC space if and only if whenever for any closed subset $A\subseteq X$, $f:(A,d)\rightarrow (\mathbb{R},|\cdot|)$ is uniformly continuous and never zero, then $\frac{1}{f}:(A,d)\rightarrow (\mathbb{R},|\cdot|)$ is uniformly continuous.
\end{theorem}
\begin{proof}
Let $(X,d)$ be a UC space, and $A$ be a closed subset of $(X,d)$. So by Proposition \ref{closed subspace of a UC space is UC}, $(A,d)$ is UC. Suppose $f:(A,d)\rightarrow (\mathbb{R},|\cdot|)$ is uniformly continuous and never zero. Then  $\frac{1}{f}:(A,d)\rightarrow (\mathbb{R},|\cdot|)$ is forward continuous on $A$. Consequently, $\frac{1}{f}:(A,d)\rightarrow (\mathbb{R},|\cdot|)$ is uniformly continuous on $A$.
			
Conversely, suppose $(X,d)$ is not a UC space. Then by Theorem \ref{UCness and uniform discreteness}, there exists a closed and discrete subset $Y \subseteq X$ which is not uniformly discrete. Thus, for every $n \in \mathbb{N}$, there exist $x_n,y_n\in Y$ such that $0 < d(x_n, y_n) < \frac{1}{n}$.  Since $Y$ is closed and discrete in $(X,d)$ and $\tau(d)=\tau(\bar{d})$, by using an argument as used in $(b)\Rightarrow (c)$ of Theorem \ref{UCness and uniform discreteness}, one can see that $(x_n)$ and $(y_n)$ do not have a cluster point in $(X,d)$. By Proposition \ref{existence of disjoint closed sets that are near}, we can assume that $x_n \neq y_m$ for all $n, m \in \mathbb{N}$. 
Take $A= \{x_{n}:n\in \mathbb{N}\}$ and $B=\{y_{n}:n\in \mathbb{N}\}$. Then $A$ and $B$ are closed subsets in $X$. Since $A\cup B$ is closed and discrete subset of $X$, for each $x_k,y_k \in A\cup B$, there exist $\delta_{x_{k}}, \delta_{y_{k}}>0$ such that $B^{d^{s}}(x_k,\delta_{\delta_{x_{k}}})\cap (A\cup B)=\{x_k\}$ and $B^{d^{s}}(y_k,\delta_{\delta_{y_{k}}})\cap (A\cup B)=\{y_k\}$, respectively. 
					
Define $f:(A\cup B,d)\rightarrow (\mathbb{R}, |\cdot|)$ as $$f(x)=\Biggl\{\begin{array}{lc}
				\frac{1}{n}, x=x_n \text{ for some } n\in \mathbb{N}, \\
				\frac{1}{n^{2}}, x=y_n \text{ for some } n\in \mathbb{N}.\end{array}$$
Clearly, $\frac{1}{f}$ is not uniformly continuous on $A\cup B$. It remains to show $f$ is uniformly continuous on $A\cup B$. Let $(a_n)$ and $(b_n)$ be sequences in $A\cup B$ such that $(a_n)$ is forward parallel to $(b_n)$, and $\epsilon >0$. Choose $n_0\in \mathbb{N}$ such that $\frac{1}{n_0}<\epsilon$. First, suppose that either both $(a_n)$ and $(b_n)$ consist of distinct terms or no term of $(a_n)$ or $(b_n)$ repeats infinitely many times.  Then for sufficiently large $n_{\epsilon}\in \mathbb{N}$, we have $a_n,b_n\in \{x_n,y_n:n\geq n_0\}$ for all $n\geq n_{\epsilon}$. By the choice of $n_0$, one can see that $|f(a_n)-f(b_n)|<\epsilon$ for all $n\geq n_{\epsilon}$. Now suppose $(a_n)$ has a term repeating infinitely many times, that is, we have a constant subsequence $a_{n_{k}}=a$ of $(a_n)$. Then $b_{n_{k}}\xrightarrow{f} a$, and consequently $b_{n_k} = a$ eventually as $B^{d^{s}}(a,\delta_{a})\cap (A\cup B)=\{a\}$.
If $a_{n_{k}}=b_{n_{k}}$ eventually, but $a_n\neq b_n$ when $n\neq n_k$, then by the argument above we again have $|f(a_n)-f(b_n)|<\epsilon$ for all $n\geq n_{\epsilon}$. Similarly, if $(b_n)$ has one term repeating infinitely many times, the same reasoning shows that $(f(a_n))\parallel (f(b_n))$. Finally, suppose $(a_n)$ has infinitely many terms repeating infinitely many times. Then if for some $1\leq i \leq n_0$, $x_i$ or $y_i$ repeats infinitely many times in $(a_n)$, then as argued above, we can find $k_i\in \mathbb{N}$ such that for all $n\geq k_i$ whenever $a_n= x_i$ or $y_i$, we have $a_n=b_n$. Let $K_{n_0}=\max\{k_i:1\leq i\leq n_0\}$. Then for all $n\geq K_{n_0}$, if $a_n=x_i$ or $a_n=y_i$ for some $1\leq i \leq n_0$, it follows that $a_n=b_n$. Similarly, since $(X,d)$ has the property backward convergence implies forward convergence, there exists $L_{n_{0}}$ such that for all $n\geq L_{n_0}$, if $b_n=x_i$ or $b_n=y_i$ for $1\leq i \leq n_0$, then $a_n=b_n$. Now set $T_{n_{0}}=\max\{ L_{n_0}, K_{n_0}\}$. Then for all $n\geq T_{n_{0}}$, we have either $a_n=b_n$ or $a_n,b_n\in \{x_n,y_n:n\geq n_0\}$. So by the earlier argument, we get $(f(a_n))$ is parallel to $(f(b_n))$.  Hence $f$ is uniformly continuous. However, $\frac{1}{f}$ is not uniformly continuous. Indeed, the sequence $(x_n)$ is forward parallel to $(y_n)$ in $A\cup B$, while the sequence $\left(\frac{1}{f(x_n)}\right)$ in not parallel to $\left(\frac{1}{f(y_n)}\right)$ .
\end{proof}	
		
The following examples justify the hypotheses of the above theorem. 		
\begin{example}
Let $(X,d)$ be a quasi-metric space as defined in Example \ref{uc implies every pseudo left k cauchy clusters}. Then $(X,d)$ is a $T_2$, UC space, but $\tau(d)\neq \tau(\bar{d})$. Now let $Y$ be a closed subset of $X$, and let $f:(X,d)\rightarrow (\mathbb{R},|\cdot|)$ be non-zero uniformly continuous map. To show $\frac{1}{f}$ is uniformly continuous, observe that if $(x_n)$ is forward parallel to $(y_n)$ then either $x_n=y_n$ eventually, or $x_n\in \{\frac{1}{n}:n\in \mathbb{N}\}$ and $y_n=0$ eventually. In the former case, the claim is straightforward. In the later case, since $f$ is uniformly continuous, we have $(f(x_n))\parallel (f(y_n))$, and hence $f(x_n)\rightarrow f(0)$. As $f$ is non-zero, it follows that $\left(\frac{1}{f(x_n)}\right)$ is parallel to $\left(\frac{1}{f(y_n)}\right)$. 
\end{example}
		
\begin{example}\label{e7}
Let $X = \mathbb{Z}_{+} \cup \{\tfrac{1}{n} : n \in \mathbb{N}\}$ and define 
$d : X \times X \to \mathbb{R}^{+} \cup \{0\}$ by 
$$
d(x,y) :=
\begin{cases}
d(x,0) + d(0,y), & \text{if } x \neq y, \\[6pt]
0, & \text{if } x = y,
\end{cases}
$$
where $d(x,0) = x$ and $d(0,y) = \tfrac{1}{y}$ for $y \neq 0$. Then $d$ is a quasi-metric on $\mathbb{R}$.
			
The space $(X,d)$ is a $T_2$ quasi-metric space. However, it does not have the property  forward convergence implies backward convergence, since the sequence $(x_n) = (n)$ is forward convergent to $0$ but not backward convergent. 
			
Define $f : (X,d) \to (\mathbb{R},|\cdot|)$ by 
$$
f(x) =
\begin{cases}
\frac{1}{x}, & \text{if } x \neq 0, \\[6pt]
0, & \text{if } x = 0.
\end{cases}
$$
			
We first show that $f$ is forward continuous. Let $(x_n)$ be any forward convergent sequence in $(X,d)$ such that $x_n \xrightarrow{f} x \in X$.  We claim that  $f(x_n) \to f(x)$. We can assume that $(x_n)$ is not eventually constant. Then only finitely many $x_n$ can belong to the set $\{\tfrac{1}{n}:n\in\mathbb{N}\}$, otherwise $(x_n)$ has to be eventually constant. So $(x_n)$ eventually lies in $\mathbb{Z}^{+}$. Since $x_n \xrightarrow{f} x$, it follows that $x=0$. So $d(0,x_n) = \frac{1}{x_n} \to 0$. Hence $f(x_n)\to 0$. However, $f$ is not uniformly continuous, since the sequence $(x_n)=(\tfrac{1}{n})$ is forward parallel to $(y_n)=0$, but $(f(x_n))$ is not parallel to $(f(y_n))$. Therefore, $(X,d)$ is not a UC space.
			
Suppose $Y\subseteq X$ is a closed subset of $X$ and $g:(Y,d)\rightarrow (\mathbb{R},|\cdot|)$ be a non-vanishing uniformly continuous function. Suppose $(x_n)$ is forward parallel to $(y_n)$ in $(Y,d)$. Then we can find $n_0 \in \mathbb{N}$ such that for any $n\geq n_0$, we have either $x_n = y_n$ or $x_n \in \{\frac{1}{n}: n \in \mathbb{N}\}$ and $y_n \in \mathbb{Z}_+$. If $x_n = y_n$ eventually, then $\left(\frac{1}{g(x_n)}\right)\parallel \left(\frac{1}{g(y_n)}\right)$. Now suppose $(x_n)$ and $(y_n)$ are not eventually equal. If $x_n=y_n$ only for finitely many $n$, then we have $x_n \in \{\frac{1}{n}: n \in \mathbb{N}\}$ and $y_n \in \mathbb{Z}_+$ eventually. Since $(x_n)$ is forward parallel to $(y_n)$, for every $\epsilon >0$, there exists $n_0\in \mathbb{N}$ such that $d(x_n,y_n)=x_n+\frac{1}{y_n}<\epsilon$ for all $n\geq n_0$. In fact $d(x_n,0)=x_n<\epsilon$ and $d(0,y_n)=\frac{1}{y_n}<\epsilon$ for all $n\geq n_0$. Hence $(x_n)$ is forward parallel to the constant sequence $z_n=0$, and the constant sequence $z_n=0$ is forward parallel to the sequence  $(y_n)$. Since $f$ is uniformly continuous, we get that $f(x_n)\rightarrow f(0)$ and $f(y_n)\rightarrow f(0)$. As $f$ is a non-zero function, $\frac{1}{f(x_n)}\rightarrow \frac{1}{f(0)}$ and  $\frac{1}{f(y_n)}\rightarrow \frac{1}{f(0)}$. Therefore, $\left(\frac{1}{f(x_n)}\right)$ is forward parallel to $\left(\frac{1}{f(y_n)}\right)$.

If $x_n=y_n$ only for infinitely many $n$, then we can split $(x_n)$ and $(y_n)$ into complementary subsequences $(x_{n_{k}})$,  $(x_{n_{m}})$,  $(y_{n_{k}})$ and  $(y_{n_{m}})$ such that $x_{n_{k}}= y_{n_{k}}$ for all $k\in \mathbb{N}$ and $x_{n_{m}}\in \{\frac{1}{n}: n \in \mathbb{N}\}$ while $y_{n_{m}}  \in \mathbb{Z}_+$ for all $m\in \mathbb{N}$. Then $\frac{1}{f(x_{n_{k}})}=\frac{1}{f(y_{n_{k}})}$ and the sequence $\left(\frac{1}{f(x_{n_{m}})}\right)$ is forward parallel to $\left(\frac{1}{f(y_{n_{m}})}\right)$. Therefore,  $\left(\frac{1}{f(x_n)}\right)$ is forward parallel to $\left(\frac{1}{f(y_n)}\right)$.
\end{example}
	
\bibliographystyle{amsplain}
\bibliography{paper2}
		
\end{document}